\documentclass{article}
\usepackage{amsmath,amsfonts,amssymb,amsthm,mathtools,bbm}
\usepackage{microtype}
\usepackage{lipsum}
\usepackage[T1]{fontenc}
\usepackage{url,hyperref}
\usepackage{indentfirst}
\hypersetup{colorlinks, linkcolor={red}, citecolor={blue}, urlcolor={black}}
\usepackage{comment}
\usepackage[margin=1.25in]{geometry}
\usepackage{enumerate,enumitem}
\usepackage{subcaption}
\usepackage{tikz}
\usetikzlibrary{calc,matrix,positioning,angles,quotes}
\usetikzlibrary{shapes.geometric}
\usetikzlibrary{patterns,decorations.pathreplacing}
\usepackage{array,multirow}

\allowdisplaybreaks

\newcommand\ceil[1]{\left\lceil#1\right\rceil}
\newcommand\floor[1]{\left\lfloor#1\right\rfloor}

\newtheorem{theorem}{Theorem}[section]
\newtheorem{lemma}[theorem]{Lemma}

\newtheorem{corollary}[theorem]{Corollary}

\newtheorem{conjecture}[theorem]{Conjecture}
\newtheorem{claim}[theorem]{Claim}

\title{Extremal persistence probabilities of\\exchangeable sign-invariant random variables}
\author{Daniel I\v{l}kovi\v{c}\thanks{Google DeepMind. Email: \url{danieli@google.com}} \and Jun Yan\thanks{Faculty of Mathematics, University of Waterloo, Canada. Email: \url{j228yan@uwaterloo.ca}.}}
\date{}

\begin{document}

\maketitle

\begin{abstract}
Let $(X_1,\ldots,X_n)$ be a random variable in $(\mathbb{R}\setminus\{0\})^n$ that is both exchangeable and sign-invariant. For every $k\in[n]$, let $S_k=\sum_{i=1}^kX_i$. Define the weak persistence probability as $\mathbb{P}(S_1,\ldots,S_n\geq0)$, and the strong persistence probability as $\mathbb{P}(S_1,\ldots,S_n>0)$.

From previous results, the optimal lower bound for $\mathbb{P}(S_1,\ldots,S_n\geq0)$ and the optimal upper bound for $\mathbb{P}(S_1,\ldots,S_n>0)$ are known. We complete the picture by determining the optimal upper bound for $\mathbb{P}(S_1,\ldots,S_n\geq0)$ and the optimal lower bound for $\mathbb{P}(S_1,\ldots,S_n>0)$ as follows.
\[\frac{1}{2^n}\binom{n-1}{\floor{(n-1)/2}}\leq\mathbb{P}(S_1,\ldots,S_n>0)\leq\frac{1}{4^n}\binom{2n}{n}\leq\mathbb{P}(S_1,\ldots,S_n\geq0)\leq\frac{1}{2^n}\binom{n}{\floor{n/2}}.\]
In particular, this implies that the weak and strong persistence probabilities of every exchangeable and sign-invariant random variable $(X_1,\ldots,X_n)$ in $(\mathbb{R}\setminus\{0\})^n$ are of the order $n^{-1/2}$. We also obtain some related results, one in the deterministic setting, and one when the random variables are allowed to take the value 0. 
\end{abstract}

\section{Introduction}\label{sec:intro}
Let $(X_1,\ldots,X_n)$ be a random variable in $\mathbb{R}^n$. Associate to $(X_1,\ldots,X_n)$ a one-dimensional random walk $W=W(X_1,\ldots,X_n)$ starting at 0, and moving by $X_k$ at time $k$ for every $k\in[n]$.
\begin{itemize}
    \item $(X_1,\ldots,X_n)$ is \textit{exchangeable} if for any permutation $\sigma\in\mathfrak{S}_n$, $(X_1,\ldots,X_n)$ has the same distribution as $(X_{\sigma(1)},\ldots,X_{\sigma(n)})$.
    \item $(X_1,\ldots,X_n)$ is \textit{sign-invariant} if for any $s\in\{-1,1\}^n$, $(X_1,\ldots,X_n)$ has the same distribution as $(s_1X_1,\ldots,s_nX_n)$.
    \item $(X_1,\ldots,X_n)$ is \textit{sum-distinct} if almost surely, for all distinct $I,J\subset [n]$, $\sum_{i\in I}X_i\not=\sum_{i\in J}X_i$. 
\end{itemize} For every $k\in[n]$, let $S_k=\sum_{i=1}^kX_i$ be the position of $W$ after time $k$. 
\begin{itemize}
    \item The \textit{weak persistence probability} of $(X_1,\ldots,X_n)$ is $\mathbb{P}(S_1,\ldots,S_n\geq0)$, the probability that $W$ stays non-negative.
    \item The \textit{strong persistence probability} of $(X_1,\ldots,X_n)$ is $\mathbb{P}(S_1,\ldots,S_n>0)$, the probability that $W$ stays positive after time 0.
\end{itemize}
We also say that $W=W(X_1,\ldots,X_n)$ is exchangeable, sign-invariant, or sum-distinct if $(X_1,\ldots,X_n)$ has the respective property, and define the weak and strong persistence probability of $W$ as the corresponding probability of $(X_1,\ldots,X_n)$.

A random variable $X$ in $\mathbb{R}$ is called \textit{atomless} if for every $x\in\mathbb{R}$, $\mathbb{P}(X=x)=0$. The following somewhat surprising result of Andersen from 1954 shows that when $X_1,\ldots,X_n$ are independent identically distributed copies of the same symmetric atomless random variable $X$, then the weak and strong persistence probabilities coincide, and are independent of the choice of $X$. 
\begin{theorem}[{\cite[Theorem 2]{A}}]\label{thm:symatom}
Let $X_1,\ldots,X_n$ be independent identically distributed copies of a symmetric atomless random variable $X$. For every $k\in[n]$, let $S_k=\sum_{i=1}^kX_i$. Then, \[\mathbb{P}(S_1,\ldots,S_n>0)=\mathbb{P}(S_1,\ldots,S_n\geq 0)=\frac{1}{4^n}\binom{2n}{n}.\] 
\end{theorem}

In fact, Andersen proved a variety of other results in~\cite{A} about partial sums of (not necessarily independent) random variables $X_1,\ldots,X_n$. For persistence probabilities, one of the most general results he proved is~{\cite[Theorem 4]{A}}, which shows that Theorem~\ref{thm:symatom} still holds under the weaker condition that $(X_1,\ldots,X_n)$ is exchangeable, sign-invariant, and sum-distinct. This result was recently rediscovered by Balister, Donderwinkel, Groenland, Johnston, and Scott in~\cite{BDGJS}, and revisited by Berger and B\'ethencourt in~\cite{BB}. They obtained the following slightly more general version.
\begin{theorem}[{\cite[Lemma 3.4]{BDGJS}}, {\cite[Theorem 2.1]{BB}}]\label{thm:esatom}
Let $(X_1,\ldots,X_n)$ be an exchangeable and sign-invariant random variable in $\mathbb{R}^n$. For every $k\in[n]$, let $S_k=\sum_{i=1}^kX_i$. Then, 
\[\mathbb{P}(S_1,\ldots,S_n>0)\leq\frac{1}{4^n}\binom{2n}{n}\leq\mathbb{P}(S_1,\ldots,S_n\geq 0).\]
Furthermore, if $(X_1,\ldots,X_n)$ is sum-distinct, then
\[\mathbb{P}(S_1,\ldots,S_n>0)=\frac{1}{4^n}\binom{2n}{n}=\mathbb{P}(S_1,\ldots,S_n\geq 0).\]
\end{theorem}

The result above establishes a tight lower bound for the weak persistence probability and a tight upper bound for the strong persistence probability. The natural question is whether one can obtain meaningful bounds in the other two directions. However, when $(X_1,\ldots,X_n)$ is exchangeable and sign-invariant, but not necessarily sum-distinct, no general meaningful upper bound for the weak persistence probability or lower bound for the strong persistence probability independent of $(X_1,\ldots,X_n)$ can exist, as when $(X_1,\ldots,X_n)$ is the random variable that is always equal to $(0,\ldots,0)$, $\mathbb{P}(S_1,\ldots,S_n\geq 0)=1$ and $\mathbb{P}(S_1,\ldots,S_n>0)=0$. 

Nevertheless, note that each step in the above pathological example is a 0-step, while in the sum-distinct setting, with probability 1 every step is non-zero. Hence, perhaps meaningful bounds can be found in terms of the random variable $M$ counting the number of $k\in[n]$ with $X_k\not=0$. In this direction, Berger and B\'ethencourt~{\cite[Remark 2.2]{BB}} proved a rudimentary lower bound $\mathbb{P}(S_1,\ldots,S_n>0)\geq\frac{\mathbb{E}(M)}{2^{2n-1}n}\binom{2n-2}{n-1}$, while Balister, Donderwinkel, Groenland, Johnston, and Scott posed the following fascinating conjecture in~\cite{BDGJS} on the upper bound of $\mathbb{P}(S_1,\ldots,S_n\geq0)$ when $M$ is always $n$, or equivalently when there is no 0-step.
\begin{conjecture}[{\cite[Conjecture 6.1]{BDGJS}}]\label{conj:BDGJS}
Let $(X_1,\ldots,X_n)$ be an exchangeable and sign-invariant random variable in $(\mathbb{R}\setminus\{0\})^n$. For every $k\in[n]$, let $S_k=\sum_{i=1}^kX_i$. Then,
\[\mathbb{P}(S_1,\ldots,S_n\geq0)\leq\frac{1}{2^n}\binom{n}{\floor{n/2}}.\]
\end{conjecture}
It can be shown (see Lemma~\ref{lemma:tight}) that when $X$ is the random variable taking values $-1$ and $1$ each with probability $1/2$, and $X_1,\ldots,X_n$ are independent identically distributed copies of $X$, then $\mathbb{P}(S_1,\ldots,S_n\geq0)=\frac{1}{2^n}\binom{n}{\floor{n/2}}$. Thus, equivalently, Conjecture~\ref{conj:BDGJS} is saying that the simple symmetric random walk maximizes the weak persistence probability under the stated assumptions. 

For the main result of this paper, we prove Conjecture~\ref{conj:BDGJS}, which establishes a tight upper bound for the weak persistence probability of an exchangeable and sign-invariant random variable $(X_1,\ldots,X_n)\in(\mathbb{R}\setminus\{0\})^n$. Moreover, we complete the picture by also obtaining a tight lower bound for the strong persistence probability in the same setting. Both bounds are achieved by the simple symmetric random walk. 
\begin{theorem}\label{thm:main1}
Let $(X_1,\ldots,X_n)$ be an exchangeable and sign-invariant random variable in $(\mathbb{R}\setminus\{0\})^n$. For every $k\in[n]$, let $S_k=\sum_{i=1}^kX_i$. Then,
\[\mathbb{P}(S_1,\ldots,S_n\geq 0)\leq\frac{1}{2^n}\binom{n}{\floor{n/2}},\quad\text{ and }\quad\mathbb{P}(S_1,\ldots,S_n>0)\geq\frac{1}{2^n}\binom{n-1}{\floor{(n-1)/2}}.\]
Furthermore, if $X_1,\ldots,X_n$ are independent identically distributed copies of the random variable $X$ taking values $-1$ and $1$ with equal probability, then both bounds are achieved.
\end{theorem}
In fact, both equality cases can be fully characterized using the full equality characterizations in Theorem~\ref{thm:mainweak} and Theorem~\ref{thm:mainstrong} below, but they are somewhat technical to state, and so are omitted. Combining Theorem~\ref{thm:esatom} and Theorem~\ref{thm:main1}, we have asymptotically that
\[\sqrt\frac{1+o(1)}{2\pi n}\leq\mathbb{P}(S_1,\ldots,S_n>0)\leq\sqrt{\frac{1+o(1)}{\pi n}}\leq\mathbb{P}(S_1,\ldots,S_n\geq 0)\leq\sqrt{\frac{2+o(1)}{\pi n}}.\]
In particular, the weak and strong persistence probabilities of every exchangeable and sign-invariant random variable in $(\mathbb{R}\setminus\{0\})^n$ are of the order $n^{-1/2}$, and are at most roughly a factor 2 apart.

To prove Theorem~\ref{thm:main1}, we follow the approach of Balister, Donderwinkel, Groenland, Johnston, and Scott~\cite{BDGJS}, and study the following closely related extremal problem. Given $x_1,\ldots,x_n>0$, define $A_\geq(x_1,\ldots,x_n)$ to be the set
\[A_\geq(x_1,\ldots,x_n)=\left\{(\sigma,s)\in\mathfrak{S}_n\times\{-1,1\}^n:\sum_{i=1}^ks_ix_{\sigma(i)}\geq0\text{ for all }k\in[n]\right\},\]
and $A_>(x_1,\ldots,x_n)$ to be the set
\[A_>(x_1,\ldots,x_n)=\left\{(\sigma,s)\in\mathfrak{S}_n\times\{-1,1\}^n:\sum_{i=1}^ks_ix_{\sigma(i)}>0\text{ for all }k\in[n]\right\}.\]
Call $(x_1,\ldots,x_n)\in\mathbb{R}^n$ \textit{sum-distinct} if for all distinct $I,J\subset[n]$, $\sum_{i\in I}x_i\not=\sum_{j\in J}x_j$. Reproducing a result of Burns in~{\cite[Section 2.3, 2.4]{B}}, Balister, Donderwinkel, Groenland, Johnston, and Scott proved in {\cite[Lemma 3.3]{BDGJS}} that
\[|A_>(x_1,\ldots,x_n)|\leq\frac{n!}{2^n}\binom{2n}{n}\leq |A_\geq(x_1,\ldots,x_n)|,\]
and that both inequalities are achieved when $(x_1,\ldots,x_n)$ is sum-distinct. Furthermore, they showed in {\cite[Lemma 3.4]{BDGJS}} that this result implies Theorem~\ref{thm:esatom}. In our case, we prove the following two results which then imply Theorem~\ref{thm:main1} by a similar argument (see Lemma~\ref{lemma:reduction}). In fact, Conjecture~\ref{conj:BDGJS} is stated using $A_\geq(x_1,\ldots,x_n)$ in~{\cite[Conjecture 6.1]{BDGJS}}.
\begin{theorem}\label{thm:mainweak}
Let $x_1,\ldots,x_n>0$. Then,
\[|A_\geq(x_1,\ldots,x_n)|\leq n!\binom{n}{\floor{n/2}}.\]
Moreover, equality is achieved, up to permutations, if and only if $x_1=\cdots=x_n$, or $n$ is odd and $x_1=\cdots=x_{n-1}=x_n/2$. 
\end{theorem}

\begin{theorem}\label{thm:mainstrong}
Let $x_1,\ldots,x_n>0$. Then,
\[|A_>(x_1,\ldots,x_n)|\geq n!\binom{n-1}{\floor{(n-1)/2}}.\]
Moreover, equality is achieved, up to permutations, if and only if $x_1=\cdots=x_n$, or $n$ is odd and $x_1=\cdots=x_{n-1}=x_n/2$. 
\end{theorem}

Surprisingly, our proofs for Theorem~\ref{thm:mainweak} and Theorem~\ref{thm:mainstrong} use entirely different methods. The proof of Theorem~\ref{thm:mainweak} is more elementary, and uses an induction assisted by decompositions of the corresponding random walks according to when they first become negative. Unfortunately, we were not able to adapt this strategy to the strong persistence setting. Instead, to prove Theorem~\ref{thm:mainstrong}, we first show that we may reduce it to the case when each of $x_1,\ldots,x_n$ is either 1 or 2, then find the minimizers among these with the help of some generating function techniques. This more complicated approach can be used to prove Theorem~\ref{thm:mainweak} as well, but we opted to present the more elegant induction proof for that instead. It would be very interesting to know if there is an easier or more elementary proof of Theorem~\ref{thm:mainstrong}.

Together, Theorem~\ref{thm:esatom} and Theorem~\ref{thm:main1} give optimal lower and upper bounds for both the weak and strong persistence probabilities of exchangeable and sign-invariant random variables $(X_1,\ldots,X_n)\in(\mathbb{R}\setminus\{0\})^n$, or equivalently when the random variable $M$ counting the number of non-zero steps is always $n$. The case when $M$ is always $0$ is trivial. With some additional arguments, we obtain the following extensions of the aforementioned theorems. The first one gives optimal lower and upper bounds for both the weak and strong persistence probabilities of exchangeable and sign-invariant random variables $(X_1,\ldots,X_n)$ in terms of $M$. The second one gives tighter bounds in the special case when $X_1,\ldots,X_n$ are independent identically distributed copies of a symmetric random variable $X$ in terms of $\mathbb{P}(X>0)$. 
\begin{corollary}\label{cor:pmain}
Let $(X_1,\ldots,X_n)$ be an exchangeable and sign-invariant random variable in $\mathbb{R}^n$. Let $M$ be the random variable counting the number of $k\in[n]$ with $X_k\not=0$. For every $k\in[n]$, let $S_k=\sum_{i=1}^kX_i$. Then,
\[\mathbb{E}\left(\frac{1}{4^M}\binom{2M}{M}\right)\leq\mathbb{P}(S_1,\ldots,S_n\geq0)\leq\mathbb{E}\left(\frac{1}{2^M}\binom{M}{\floor{M/2}}\right),\]
\[\frac1n\mathbb{E}\left(\frac{M}{2^M}\binom{M-1}{\floor{(M-1)/2}}\right)\leq\mathbb{P}(S_1,\ldots,S_n>0)\leq\frac1n\mathbb{E}\left(\frac{M}{4^M}\binom{2M}{M}\right),\]
where we interpret $\frac{M}{2^M}\binom{M-1}{\floor{(M-1)/2}}$ as 0 when $M$ takes the value 0. Moreover, for any given random variable $M$, all bounds can be achieved.
\end{corollary}

\begin{corollary}\label{cor:piid}
Let $X$ be a symmetric random variable with $\mathbb{P}(X>0)=p$. Let $X_1,\ldots,X_n$ be independent identically distributed copies of $X$. For every $k\in[n]$, let $S_k=\sum_{i=1}^kX_i$. Then,
\begin{align*}
\mathbb{P}(S_1,\ldots,S_n\geq0)&\geq\sum_{m=0}^n\binom nm(2p)^m(1-2p)^{n-m}\frac1{4^m}\binom{2m}{m},\\
\mathbb{P}(S_1,\ldots,S_n\geq0)&\leq\sum_{m=0}^n\binom nm(2p)^m(1-2p)^{n-m}\frac{1}{2^m}\binom{m}{\floor{m/2}},\\
\mathbb{P}(S_1,\ldots,S_n>0)&\geq\sum_{m=1}^n\binom{n-1}{m-1}(2p)^m(1-2p)^{n-m}\frac1{2^m}\binom{m-1}{\floor{(m-1)/2}},\\
\mathbb{P}(S_1,\ldots,S_n>0)&\leq\sum_{m=1}^n\binom{n-1}{m-1}(2p)^m(1-2p)^{n-m}\frac{1}{4^m}\binom{2m}{m}.
\end{align*}
Moreover, for every $0\leq p\leq 1/2$, all bounds can be achieved.
\end{corollary}

The rest of the paper is organized as follows. In Section~\ref{sec:prelim}, we prove two preliminary results. The first one verifies the equality cases claimed in Theorem~\ref{thm:mainweak} and Theorem~\ref{thm:mainstrong}, and the second one shows that these two theorems imply our main result Theorem~\ref{thm:main1}. We then prove Theorem~\ref{thm:mainweak} and Theorem~\ref{thm:mainstrong} in Section~\ref{sec:weak} and Section~\ref{sec:strong}, respectively. Finally, we prove Corollary~\ref{cor:pmain} and Corollary~\ref{cor:piid} in Section~\ref{sec:cor}. 

\section{Preliminary results}\label{sec:prelim}
\begin{lemma}\label{lemma:tight}
If $x_1=\cdots=x_n>0$, or $n$ is odd and $x_1=\cdots=x_{n-1}=x_n/2>0$, then 
\[|A_\geq(x_1,\ldots,x_n)|=n!\binom{n}{\floor{n/2}},\quad\text{ and }\quad|A_>(x_1,\ldots,x_n)|=n!\binom{n-1}{\floor{(n-1)/2}}.\]
\end{lemma}
\begin{proof}
First suppose without loss of generality that $x_1=\cdots=x_n=1$. Let $\mathbbm{1}_n$ be the length-$n$ vector of $1$'s. Let $W_n$ be the set of $n$-step walks starting at 0, with each step being $-1$ or 1, and never going below $0$. Then, $|A_\geq(\mathbbm{1}_n)|=n!|W_n|$. From the well-known ballot problem, the number of such walks with $p$ 1-steps is $\binom{n}{p}-\binom{n}{p+1}$. Summing over $p$ from $\ceil{n/2}$ to $n$ gives $|A_\geq(\mathbbm{1}_n)|=n!|W_n|=n!\binom{n}{\ceil{n/2}}=n!\binom{n}{\floor{n/2}}$. For an $n$-step walk in $W_n$ that stays strictly positive, the first step must be a 1-step, and the remaining $(n-1)$ steps form a walk in $W_{n-1}$. Thus, $|A_>(\mathbbm{1}_n)|=n!|W_{n-1}|=n!\binom{n-1}{\floor{(n-1)/2}}$.

Now suppose that $n$ is odd, and without loss of generality that $x_1=\cdots=x_{n-1}=x_n/2=1$. We will only prove $|A_\geq(1,\ldots,1,2)|=n!\binom{n}{\floor{n/2}}$ since the proof for $|A_>(1,\ldots,1,2)|$ is very similar. Let $W_n'$ be the set of $n$-step walks starting at 0, with $n-1$ steps being $-1$ or 1, the other step being $-2$ or 2, and never going below $0$. Then, $|A_\geq(1,\ldots,1,2)|=(n-1)!|W_n'|$.

Let $MW_{n+1}$ be the set of walks in $W_{n+1}$ with one instance of two consecutive $(-1)$-steps or two consecutive 1-steps marked. Note that there is a bijection between $W_n'$ and $MW_{n+1}$, since we can expand the $(-2)$-step or the $2$-step into two consecutive marked $(-1)$-steps or 1-steps in one direction, and combine the two marked steps into one in the other direction. Thus, it suffices to show that $|MW_{n+1}|=n\binom{n}{\floor{n/2}}$.

For a walk $w\in W_{n+1}$, let $c(w)$ be the number of instances of two consecutive $(-1)$-steps or two consecutive $1$-steps. Let $p(w)$ be the number of peaks in $w$, that is the number of instances of a $1$-step followed immediately by a $(-1)$-step. Let $s(w)$ be the number of sign changes in $w$. Then, $s(w)=2p(w)$ if $w$ ends with a 1-step, and $s(w)=2p(w)-1$ if $w$ ends with a $(-1)$-step. Note that exactly half of the walks in $W_{n+1}$ end with a 1-step, since $n$ is odd implies that after $n$ steps the walk cannot be at 0, and thus both choices are possible for the last step. Also, $c(w)=n-s(w)$ because each of the $n$ pairs of consecutive steps of $w$ either has the same sign or has a sign change. Putting these together, we get 
\[|MW_{n+1}|=\sum_{w\in W_{n+1}}c(w)=n|W_{n+1}|-\sum_{w\in W_{n+1}}s(w)=\frac{2n+1}2|W_{n+1}|-2\sum_{w\in W_{n+1}}p(w).\]

To finish, we claim that $\sum_{w\in W_{n+1}}p(w)=n|W_{n-1}|$. Indeed, take any $w\in W_{n+1}$, we can remove one of its peaks to get $w'\in W_{n-1}$, with the original peak location marked, and given any $w'\in W_{n-1}$ with one of its $n$ gaps between steps (including the beginning and end) marked, we can insert a peak there to obtain $w\in W_{n+1}$. This gives a bijection proving the claim. Therefore, using from above that $|W_n|=\binom{n}{\floor{n/2}}$ and $n$ is odd, we have
\[|MW_{n+1}|=\frac{2n+1}2|W_{n+1}|-2n|W_{n-1}|=\frac{2n+1}2\binom{n+1}{(n+1)/2}-2n\binom{n-1}{(n-1)/2}=n\binom{n}{\floor{n/2}},\]
as required.
\end{proof}

\begin{lemma}\label{lemma:reduction}
Theorem~\ref{thm:mainweak} and Theorem~\ref{thm:mainstrong} imply Theorem~\ref{thm:main1}.
\end{lemma}
\begin{proof}
Since $(X_1,\ldots,X_n)$ is exchangeable and sign-invariant, if $\sigma\in\mathfrak{S}_n$ and $s\in\{-1,1\}^n$ are both chosen uniformly at random, then with $S_k^{\sigma,s}=\sum_{i=1}^ks_i|X_{\sigma(i)}|$ for every $k\in[n]$, we have
\[\mathbb{P}(S_1,\ldots,S_n\geq0)=\mathbb{P}(S_1^{\sigma,s},\ldots,S_n^{\sigma,s}\geq0),\quad\text{ and }\quad\mathbb{P}(S_1,\ldots,S_n>0)=\mathbb{P}(S_1^{\sigma,s},\ldots,S_n^{\sigma,s}>0).\]
For every realization of $(X_1,\ldots,X_n)=(x_1,\ldots,x_n)$, by Theorem~\ref{thm:mainweak} and Theorem~\ref{thm:mainstrong},
\[\mathbb{P}(S_1^{\sigma,s},\ldots,S_n^{\sigma,s}\geq0\mid X_1=x_1,\ldots,X_n=x_n)=\frac1{n!2^n}|A_\geq(|x_1|,\ldots,|x_n|)|\leq\frac1{2^n}\binom{n}{\floor{n/2}},\]
\[\mathbb{P}(S_1^{\sigma,s},\ldots,S_n^{\sigma,s}>0\mid X_1=x_1,\ldots,X_n=x_n)=\frac1{n!2^n}|A_>(|x_1|,\ldots,|x_n|)|\geq\frac1{2^n}\binom{n-1}{\floor{(n-1)/2}}.\]
As these hold for every realization, the result follows. 
\end{proof}

\section{Proof of Theorem~\ref{thm:mainweak}}\label{sec:weak}
In this section, we prove Theorem~\ref{thm:mainweak}. We begin by defining most of the notations used in the proof, and giving a proof sketch. 

Let $x_1,\ldots,x_n>0$, let $\vec{x}=(x_1,\ldots,x_n)$, and let $(\sigma,s)\in\mathfrak{S}_n\times\{-1,1\}^n$. Then, the set considered in Theorem~\ref{thm:mainweak} can be interpreted as a subset of the $n$-step one-dimensional walks, whose step sizes are obtained from a signed permutation of $x_1,\ldots,x_n$.
\begin{itemize}
    \item Let $W=W^{\sigma,s}(\vec{x})$ be the one-dimensional walk starting at 0, and taking a step $s_kx_{\sigma(k)}$ at time $k$ for every $k\in[n]$. 
    \item For every $k\in[n]$, let $S_k^{\sigma,s}(\vec{x})=\sum_{i=1}^ks_ix_{\sigma(i)}$ be the position after time $k$ of the walk $W^{\sigma,s}(\vec{x})$.
\end{itemize}
In this language, the elements in $A_\geq(\vec{x})$ correspond to walks that stay non-negative.
\begin{itemize}
    \item Let $A_\geq(\vec{x})=\{(\sigma,s)\in\mathfrak{S}_n\times\{-1,1\}^n:S^{\sigma,s}_k(\vec{x})\geq0\text{ for all }k\in[n]\}$.
\end{itemize}
Equivalently, $A_\geq(\vec{x})$ corresponds to walks with a (not necessarily strict) global minimum at time 0.

Next, consider $A'_\geq(\vec{x})$, whose elements correspond to walks with a global minimum at time $n$.
\begin{itemize}
    \item Let $A'_\geq(\vec{x})=\{(\sigma,s)\in\mathfrak{S}_n\times\{-1,1\}^n:S^{\sigma,s}_k(\vec{x})\geq S^{\sigma,s}_n(\vec{x})\text{ for all }0\leq k\leq n\}$.
\end{itemize}
Note that by reversing the permutation $\sigma$, or equivalently reversing the walk, we get a bijection between $A_\geq(\vec{x})$ and $A'_\geq(\vec{x})$, so $|A_\geq(\vec{x})|=|A'_\geq(\vec{x})|$.

We now take the complement of the union of the sets $A_\geq(\vec{x})$ and $A'_\geq(\vec{x})$ considered above. This new set turns out to be the key to our proof. 
\begin{itemize}
    \item Let $I_\geq(\vec{x})=(\mathfrak{S}_n\times\{-1,1\}^n)\setminus(A_\geq(\vec{x})\cup A'_\geq(\vec{x}))$.
\end{itemize}
Equivalently, elements in $I_\geq(\vec{x})$ correspond to the walks whose global minimums all occur between time 1 and time $n-1$.

We are now ready to give a sketch of the proof that $|A_\geq(\vec{x})|\leq n!\binom{n}{\floor{n/2}}$. Note that by Lemma~\ref{lemma:tight}, $|A_\geq(\mathbbm{1}_n)|=n!\binom{n}{\floor{n/2}}$, so it suffices to prove that $|A_\geq(\vec{x})|\leq|A_\geq(\mathbbm{1}_n)|$. We use induction on $n$ to prove the stronger statement that $|A_\geq(\vec{x})|\leq|A_\geq(\mathbbm{1}_n)|$ and $|I_\geq(\vec{x})|\geq|I_\geq(\mathbbm{1}_n)|$. 

Assuming both of these hold for shorter vectors, we first show in Claim~\ref{claim:minI} that $|I_\geq(\vec{x})|\geq|I_\geq(\mathbbm{1}_n)|$. Let $(\sigma,s)\in I_\geq(\vec{x})$, and let $W=W^{\sigma,s}(\vec{x})$. Since $(\sigma,s)$ is not in $A_\geq(\vec{x})$, there exists a smallest $k\in[n]$ such that $W$ is negative after time $k$. Since $(\sigma,s)$ is not in $A'_\geq(\vec{x})$, the part of the walk in $W$ after time $k$ is not a walk with a global minimum at the end, or equivalently not a walk corresponding to $A'_\geq(\vec{y})$, where $\vec{y}$ consists of the step sizes in $\vec{x}$ not used in the first $k$ steps. Recall that $|A'_\geq(\vec{y})|=|A_\geq(\vec{y})|$ via a reversal operation, so the walk decomposition described above according to the first negative time allows us to write $|I_\geq(\vec{x})|$ in terms of quantities of the form $|A_\geq(\vec{y})|$ with $\vec{y}$ having shorter lengths. We can then apply induction hypothesis to bound $|A_\geq(\vec{y})|$ by $|A_\geq(\mathbbm{1}_{|\vec{y}|})|$, and finish the proof of Claim~\ref{claim:minI} with a summation by parts argument. 

To complete the proof that $|A_\geq(\vec{x})|\leq|A_\geq(\mathbbm{1}_n)|$, we consider the last step of the walk $W$. Suppose $x_k$ is the only unused step size. Then, the last step can be either $x_k$ or $-x_k$, unless after time $n-1$ the walk $W$ ends below height $x_k$. Since $x_k>0$, the number of such walks is at least the number of non-negative walks that end at 0, whose count can be represented using a simple inclusion-exclusion argument in terms of $|I_\geq(\vec{y})|$, where $\vec{y}$ consists of all step sizes in $\vec{x}$ except $x_k$. By induction hypothesis, $|I_\geq(\vec{y})|\geq|I_\geq(\mathbbm{1}_{n-1})|$, which then allows us to finish the proof. 

We now follow this sketch and prove Theorem~\ref{thm:mainweak}.
\begin{proof}[Proof of Theorem~\ref{thm:mainweak}]
We use induction to prove that for every $n\geq1$ and all $x_1,\ldots,x_n>0$, $|A_\geq(\vec{x})|\leq|A_\geq(\mathbbm{1}_n)|=n!\binom{n}{\floor{n/2}}$ with equality cases exactly as stated, and additionally $|I_\geq(\vec{x})|\geq|I_\geq(\mathbbm{1}_n)|$. The base cases when $n=1,2,3$ are easy to verify, so we assume that $n\geq 4$.
We first prove that $|I_\geq(\vec{x})| \geq |I_\geq(\mathbbm{1}_n)|$.
\begin{claim}\label{claim:minI}
$|I_\geq(\vec{x})| \geq |I_\geq(\mathbbm{1}_n)|$.
\end{claim}
\begin{proof}
For every $k\in[n]$, let $I_\geq^k(\vec{x})$ be the set of $(\sigma,s)$ such that the corresponding walk $W_{\sigma,s}(\vec{x})$ first becomes negative at time $k$. Then, 
\[(\mathfrak{S}_n\times\{-1,1\}^n)\setminus A_\geq(\vec{x})=\bigsqcup_{k=1}^n I_\geq^k(\vec{x}).\]  For each $J\subset[n]$, let $\vec{x}_J=(x_j)_{j\in J}$, and let $F_\geq(\vec{x}_J)$ be the set of $(\sigma,s)\in\mathfrak{S}_{|J|}\times\{-1,1\}^{|J|}$ such that the corresponding $|J|$-step walk first becomes negative at time $|J|$. Then, we have
\[|I_\geq^k(\vec{x})|=\sum_{J\in\binom{[n]}{k}}|F_\geq(\vec{x}_J)|(n-k)!2^{n-k},\]
since each walk counted by $|I_\geq^k(\vec{x})|$ must start with an initial $k$-step walk counted by $|F_\geq(\vec{x}_J)|$ for some $J\subset[n]$ of size $k$, then there is no restriction on the remaining $n-k$ steps. 

Similarly, using $|A'_\geq(\vec{y})|=|A_\geq(\vec{y})|$ for every $\vec{y}$ by reversal, we have 
\[|I_\geq^k(\vec{x})\cap A'_\geq(\vec{x})|=\sum_{J\in\binom{[n]}{k}} |F_\geq(\vec{x}_J)| |A'_\geq(\vec{x}_{[n]\setminus J})| = \sum_{J\in\binom{[n]}{k}} |F_\geq(\vec{x}_J)||A_\geq(\vec{x}_{[n]\setminus J})|,\]
because for any walk counted by $I_\geq^k(\vec{x})\cap A'_\geq(\vec{x})$, the only restriction on its terminal $(n-k)$-step walk is that it must have its global minimum achieved at the end, which is counted by $|A'_\geq(\vec{x}_{[n]\setminus J})|$.

Let $A^c(\vec{y})=(\mathfrak{S}_{|\vec{y}|}\times\{-1,1\}^{|\vec{y}|})\setminus A_\geq(\vec{y})$ for every $\vec{y}$. Combining the two equations above, we get
\begin{align*}
|I_\geq(\vec{x})|&=|(\mathfrak{S}_n\times\{-1,1\}^n)\setminus(A_\geq(\vec{x})\cup A'_\geq(\vec{x}))|\\
&=|A_\geq^c(\vec{x})|-|A_\geq^c(\vec{x})\cap A'_\geq(\vec{x})|\\
&=\sum_{k=1}^n\left(|I^k_\geq(\vec{x})|-|I_\geq^k(\vec{x})\cap A'_\geq(\vec{x})|\right)=\sum_{\varnothing\not=J\subset[n]}|F_\geq(\vec{x}_J)||A^c_\geq(\vec{x}_{[n]\setminus J})|.
\end{align*}

For every $0\leq k\leq n-1$, let $p^c_\geq(k)=|A^c_\geq(\mathbbm{1}_k)|/k!2^k$, noting that $p^c_\geq(0)=0$. By Lemma~\ref{lemma:tight}, $p^c_\geq(k)$ is equal to $1-\binom{k}{\floor{k/2}}/2^k$ for every $k\in[n-1]$, which is non-decreasing in $k$. For every non-empty $J\subset[n]$, we have $|A^c_\geq(\vec{x}_{[n]\setminus J})|\geq|A^c_\geq(\mathbbm{1}_{n-|J|})|$ by induction hypothesis. Therefore, by grouping according to the size of $J$, we get
\begin{equation}\label{eq:tight}
|I_\geq(\vec{x})|\geq\sum_{\varnothing\not=J\subset[n]}|F_\geq(\vec{x}_J)||A^c_\geq(\mathbbm{1}_{n-|J|})|=\sum_{k=1}^n \frac{|I_\geq^k(\vec{x})|}{(n-k)! 2^{n-k}}|A^c_\geq(\mathbbm{1}_{n-k})|=\sum_{k=1}^n|I_\geq^k(\vec{x})|p^c_\geq(n-k).
\end{equation}
Note that when $\vec{x}=\mathbbm{1}_n$, the inequality above is actually tight. 

For every $0\leq k\leq n$, let $I_\geq^{\leq k}(\vec{x})=\cup_{i=1}^k I_\geq^i(\vec{x})$, and let $I_\geq^{>k}(\vec{x})=(\mathfrak{S}_n\times\{-1,1\}^n)\setminus I_\geq^{\leq k}(\vec{x})$. Then, $I_\geq^{>k}(\vec{x})$ corresponds exactly to walks that remain non-negative for the first $k$ steps, and so
\[|I_\geq^{>k}(\vec{x})|=\sum_{J\in\binom{[n]}{k}}|A_\geq(\vec{x}_J)|(n-k)!2^{n-k}.\] 
If $k<n$, then by induction hypothesis, 
\[|I_\geq^{>k}(\vec{x})|=\sum_{J\in\binom{[n]}{k}}|A_\geq(\vec{x}_J)|(n-k)!2^{n-k}\leq\sum_{J\in\binom{[n]}{k}}|A_\geq(\mathbbm{1}_k)|(n-k)!2^{n-k}=|I_\geq^{>k}(\mathbbm{1}_n)|.\] 
It then follows by taking complement that $|I_\geq^{\leq k}(\vec{x})|\geq|I_\geq^{\leq k}(\mathbbm{1}_n)|$.

Putting everything together, and using summation by parts, we get
\begin{align*}
|I_\geq(\vec{x})|&\geq\sum_{k=1}^n|I_\geq^k(\vec{x})|p^c_\geq(n-k)\\
&=\sum_{k=1}^n\left(|I_\geq^{\leq k}(\vec{x})|-|I_\geq^{\leq k-1}(\vec{x})|\right)p^c_\geq(n-k)\\
&=-|I_\geq^{\leq 0}(\vec{x})|p^c_\geq(n-1)+\sum_{k=1}^{n-1}|I_\geq^{\leq k}(\vec{x})|\left(p^c_\geq(n-k)-p^c_\geq(n-k-1)\right) + |I_\geq^{\leq n}(\vec{x})|p^c_\geq(0)\\
&=\sum_{k=1}^{n-1}|I_\geq^{\leq k}(\vec{x})|\left(p^c_\geq(n-k)-p^c_\geq(n-k-1)\right)\\
&\geq\sum_{k=1}^{n-1}|I_\geq^{\leq k}(\mathbbm{1}_n)|\left(p^c_\geq(n-k)-p^c_\geq(n-k-1)\right)\\
&=\sum_{k=1}^n|I_\geq^k(\mathbbm{1}_n)|p^c_\geq(n-k)=|I_\geq(\mathbbm{1}_n)|,
\end{align*}
where in the last inequality we used that $p^c_\geq(k)$ is non-decreasing in $k$. This proves the claim. 
\renewcommand{\qedsymbol}{$\boxdot$}
\end{proof}
\renewcommand{\qedsymbol}{$\square$}

Let $Z_\ge(\vec{x})=A_\geq(\vec{x})\cap A'_\geq(\vec{x})$, and note that $Z_\ge(\vec{x})$ is exactly the set of $(\sigma,s)$ in $A_\geq(\vec{x})$ with $S_n(\sigma,s)=0$. Equivalently, the walks corresponding to elements in $Z_\geq(\vec{x})$ are those that start and end at 0, and stay non-negative throughout. By inclusion-exclusion, and using $|A_\geq(\vec{x})|=|A'_\geq(\vec{x})|$,
\begin{align*}
|I_\geq(\vec{x})|=|(\mathfrak{S}_n\times\{-1,1\}^n)\setminus(A_\geq(\vec{x})\cup A'_\geq(\vec{x}))|&=n!2^n-|A_\geq(\vec{x})\cup A'_\geq(\vec{x})|\\
&=n!2^n-|A_\geq(\vec{x})|-|A'_\geq(\vec{x})|+|Z_\ge(\vec{x})|\\
&=n!2^n-2|A_\geq(\vec{x})|+|Z_\ge(\vec{x})|.
\end{align*}
 
For every $i\in[n]$, let 
\[A^i_\geq(\vec{x})=\{(\sigma,s)\in A_\geq(\vec{x}_{[n]\setminus\{i\}}):S^{\sigma,s}_{n-1}(\vec{x}_{[n]\setminus\{i\}})<x_i\}.\]
In other words, the elements in this set correspond to $(n-1)$-step walks using all step sizes $x_1,\ldots,x_n$ except $x_i$, which are non-negative throughout, and end at a height less than $x_i$. Since this includes walks ending at 0, we have that in particular,  
\[|A^i_\geq(\vec{x})|\ge|Z_\ge(\vec{x}_{[n]\setminus\{i\}})|.\]

Observe that every walk corresponding to an element in $A^i_\geq(\vec{x})$ can only be extended to a walk corresponding to an element in $A_\geq(\vec{x})$ by taking an $x_i$-step, while every walk corresponding to an element in $A_\geq(\vec{x}_{[n]\setminus\{i\}})\setminus A^i_\geq(\vec{x})$ can be extended in two ways, with either an $x_i$-step or a $(-x_i)$-step. Therefore, by splitting the count of $|A_\geq(\vec{x})|$ according to the value of $\sigma(n)$, we have from induction hypothesis and Claim~\ref{claim:minI} that
\begin{align*}
|A_\geq(\vec{x})|&=\sum_{i=1}^n\left(|A^i_\geq(\vec{x})|+2|A_\geq(\vec{x}_{[n]\setminus\{i\}})|-2|A^i_\geq(\vec{x})|\right)\\
&=\sum_{i=1}^n\left(2|A_\geq(\vec{x}_{[n]\setminus\{i\}})|-|A^i_\geq(\vec{x})|\right)\\
&\leq\sum_{i=1}^n\left(2|A_\geq(\vec{x}_{[n]\setminus\{i\}})|-|Z_\ge(\vec{x}_{[n]\setminus\{i\}})|\right)\\
&=\sum_{i=1}^n\left((n-1)!2^{n-1}-|I_\ge(\vec{x}_{[n]\setminus\{i\}})|\right)\\
&\leq\sum_{i=1}^n\left((n-1)!2^{n-1}-|I_\ge(\mathbbm{1}_{n-1})|\right)=|A_\geq(\mathbbm{1}_n)|,
\end{align*}
where the last equality follows as both inequalities in this derivation are actually equalities for $\vec{x}=\mathbbm{1}_n$. Since $|A_\geq(\mathbbm{1}_n)|=n!\binom{n}{\floor{n/2}}$, this completes the proof that $|A_\geq(\vec{x})|\leq n!\binom{n}{\floor{n/2}}$.

Now we determine precisely the set of $\vec{x}\in(\mathbb{R}_{>0})^n$ such that $|A_\geq(\vec{x})|=n!\binom{n}{\floor{n/2}}$. From the last inequality above, for this to hold, we need $|I_\ge(\vec{x}_{[n]\setminus\{i\}})|=|I_\ge(\mathbbm{1}_{n-1})|$ for every $i\in[n]$. Tracing back to~(\ref{eq:tight}) in the proof of Claim~\ref{claim:minI} and taking $J=\{j\}$, we need in particular that for all distinct $i,j\in[n]$, $|A_\geq(\vec{x}_{[n]\setminus\{i,j\}})|=|A_\geq(\mathbbm{1}_{n-2})|$. 

If $n$ is even, and there are at least two distinct values among $x_1,\ldots,x_n$, then using $n\geq4$, there exist distinct $i,j\in[n]$ such that $\vec{x}_{[n]\setminus\{i,j\}}$ contains at least two distinct values. By the induction hypothesis, this implies that $|A_\geq(\vec{x}_{[n]\setminus\{i,j\}})|<|A_\geq(\mathbbm{1}_{n-2})|$, a contradiction. Therefore, by Lemma~\ref{lemma:tight}, $|A_\geq(\vec{x})|=|A_\geq(\mathbbm{1}_n)|$ if and only if $x_1=\cdots=x_n$. 

For the case when $n$ is odd, first suppose that there are at least three distinct values among $x_1,\ldots,x_n$. Using $n\geq5$, there exist $i,j\in[n]$ such that $\vec{x}_{[n]\setminus\{i,j\}}$ contains at least three distinct values, so by the induction hypothesis, $|A_\geq(\vec{x}_{[n]\setminus\{i,j\}})|<|A_\geq(\mathbbm{1}_{n-2})|$, a contradiction. Now suppose there are exactly two distinct values among $x_1,\ldots,x_n$, and the larger value appears at least twice. Then, since $n\geq 5$ is odd, one of the values appears at least three times, so there exist distinct $i,j\in[n]$ such that $\vec{x}_{[n]\setminus\{i,j\}}$ still contain two distinct values, and the larger value appears at least twice. Again, this implies that $|A_\geq(\vec{x}_{[n]\setminus\{i,j\}})|<|A_\geq(\mathbbm{1}_{n-2})|$ by the induction hypothesis, a contradiction. If there are exactly two distinct values among $x_1,\ldots,x_n$, and the larger value appears exactly once, but the ratio between the two values is not $2$, then we again get a contradiction using $n\geq5$ and the induction hypothesis by removing two copies of the smaller value. Therefore, by Lemma~\ref{lemma:tight}, $|A_\geq(\vec{x})|=|A_\geq(\mathbbm{1}_n)|$ if and only if $\vec{x}$ satisfies $x_1=\cdots=x_n$ or $x_1=\cdots=x_{n-1}=x_n/2$, up to a suitable index permutation.
\end{proof}

\section{Proof of Theorem~\ref{thm:mainstrong}}\label{sec:strong}
In this section, we prove Theorem~\ref{thm:mainstrong}. The proof consists of three main steps which will be carried out in separate subsections. First, in Section~\ref{sec:reduce}, we reduce the problem to the case when $\vec{x}=\mathbbm{1}_p\mathbbm{2}_q$ is a vector of 1's and 2's. Then, in Section~\ref{sec:exact}, we use generating function techniques to obtain an exact formula for $|A_>(\mathbbm{1}_p\mathbbm{2}_q)|$. Finally, in Section~\ref{sec:opt}, we determine for fixed $n=p+q$, which $p$ and $q$ minimize $|A_>(\mathbbm{1}_p\mathbbm{2}_q)|$.

The crucial starting point in our proof is the following special case of what is often referred to as the Bohnenblust-Spitzer identity. In what follows, given a permutation $\sigma\in\mathfrak{S}_n$, let $C(\sigma)$ be the set of subsets of $[n]$ that form a cycle in the cycle decomposition of $\sigma$. In particular, the sets contained in $C(\sigma)$ form a partition of $[n]$. 
\begin{lemma}[{\cite[Theorem 2.3]{S}}]\label{lemma:spitzer}
Let $x_1,\ldots,x_n\in\mathbb{R}$. Then,
\[\sum_{\sigma\in\mathfrak{S}_n}\prod_{k=1}^n\mathbbm{1}\left(\sum_{i=1}^kx_{\sigma(i)}>0\right)=\sum_{\sigma\in\mathfrak{S}_n}\prod_{C\in C(\sigma)}\mathbbm{1}\left(\sum_{i\in C}x_i>0\right).\]
\end{lemma}

Given $x_1,\ldots,x_n>0$, for each $C\subset[n]$, let $E_C(x_1,\ldots,x_n)$ be the set of $s\in\{-1,1\}^C$ satisfying $\sum_{i\in C}s_ix_i=0$. When $C=[n]$, we abbreviate $E_{[n]}(x_1,\ldots,x_n)$ to $E(x_1,\ldots,x_n)$. A straightforward corollary of Lemma~\ref{lemma:spitzer} for our purpose is the following. 
\begin{corollary}\label{cor:spitzer}
Let $x_1,\ldots,x_n>0$. Then, 
\[|A_>(x_1,\ldots,x_n)|=\sum_{\sigma\in\mathfrak{S}_n}\prod_{C\in C(\sigma)}\frac{2^{|C|}-|E_C(x_1,\ldots,x_n)|}2.\]
\end{corollary}
\begin{proof}
For every $C\subset[n]$ and every $s\in\{-1,1\}^C$, observe that by symmetry,
\[\sum_{s\in\{-1,1\}^C}\mathbbm{1}\left(\sum_{i\in C}s_ix_i>0\right)=\frac{2^{|C|}-|E_C(x_1,\ldots,x_n)|}2.\]
By Lemma~\ref{lemma:spitzer}, for every $s\in\{-1,1\}^n$, we have
\[\sum_{\sigma\in\mathfrak{S}_n}\prod_{k=1}^n\mathbbm{1}\left(\sum_{i=1}^ks_{\sigma(i)}x_{\sigma(i)}>0\right)=\sum_{\sigma\in\mathfrak{S}_n}\prod_{C\in C(\sigma)}\mathbbm{1}\left(\sum_{i\in C}s_ix_i>0\right).\]
Putting these together, we get
\begin{align*}
|A_>(x_1,\ldots,x_n)|&=\sum_{s\in\{-1,1\}^n}\sum_{\sigma\in\mathfrak{S}_n}\prod_{k=1}^n\mathbbm{1}\left(\sum_{i=1}^ks_ix_{\sigma(i)}>0\right)\\
&=\sum_{s\in\{-1,1\}^n}\sum_{\sigma\in\mathfrak{S}_n}\prod_{k=1}^n\mathbbm{1}\left(\sum_{i=1}^ks_{\sigma(i)}x_{\sigma(i)}>0\right)\\
&=\sum_{\sigma\in\mathfrak{S}_n}\sum_{s\in\{-1,1\}^n}\prod_{C\in C(\sigma)}\mathbbm{1}\left(\sum_{i\in C}s_ix_i>0\right)\\
&=\sum_{\sigma\in\mathfrak{S}_n}\prod_{C\in C(\sigma)}\left(\sum_{s\in\{-1,1\}^C}\mathbbm{1}\left(\sum_{i\in C}s_ix_i>0\right)\right)\\
&=\sum_{\sigma\in\mathfrak{S}_n}\prod_{C\in C(\sigma)}\frac{2^{|C|}-|E_C(x_1,\ldots,x_n)|}2,
\end{align*}
as required.
\end{proof}

While Lemma~\ref{lemma:spitzer} was originally stated and proved combinatorially, we take a small digression here to remark that this is in fact only a special case of a far more general phenomenon in algebra. Let $K$ be a field of characteristic 0, let $\lambda\in K$, and let $A$ be a commutative $K$-algebra. Originating from the work of Baxter, a linear operator $R$ on $A$ is called a \textit{Rota-Baxter operator} of \textit{weight} $\lambda$ if for all $x,y\in A$,
\[R(x)R(y)=R(R(x)y)+R(xR(y))+\lambda R(xy).\]
Rota and Smith~\cite{RS} showed in 1972 that for any $x_1,\ldots,x_n\in A$, 
\begin{equation}\label{eq:rs}
\sum_{\sigma\in\mathfrak{S}_n}\prod_{k=1}^nR\left(\prod_{i=1}^kx_{\sigma(i)}\right)=\sum_{\sigma\in\mathfrak{S}_n}(-\lambda)^{n-|C(\sigma)|}\prod_{C\in C(\sigma)}R\left(\prod_{i\in C}x_i\right).
\end{equation}
Furthermore, non-commutative analogues of this result have also been obtained by Ebrahimi-Fard, Manchon, and Patras~\cite{EMP}.

In the language of Rota-Baxter operators, Lemma~\ref{lemma:spitzer} can be recovered as follows. Let $G$ be the additive group whose elements are integer combinations of $x_1,\ldots,x_n$. Let $A=K[G]$ be the group algebra with basis $\{e_g:g\in G\}$. Let $R$ be the linear operator defined by $R(e_g)=e_g$ if $g>0$, and $R(e_g)=0$ otherwise. Then, it is straightforward to check that $R$ is a Rota-Baxter operator with weight $-1$, and Lemma~\ref{lemma:spitzer} follows from~(\ref{eq:rs}).

\subsection{Reduction}\label{sec:reduce}
In this subsection, we reduce Theorem~\ref{thm:mainstrong} to the case when each of $x_1,\ldots,x_n$ is either 1 or 2. First, we use basic linear algebra to reduce to the case when each of $x_1,\ldots,x_n$ is an integer. 
\begin{lemma}\label{lemma:reduceint}
For any $x_1,\ldots,x_n>0$, there exist positive integers $y_1,\ldots,y_n$ such that for every $C\subset[n]$,
$E_C(x_1,\ldots,x_n)=E_C(y_1,\ldots,y_n)$. In particular, $|A_>(x_1,\ldots,x_n)|=|A_>(y_1,\ldots,y_n)|$.
\end{lemma}
\begin{proof}
For variables $v_1,\ldots,v_n\in\mathbb{R}$, let $H_1$ be the set of all hyperplanes of the form $\sum_{i\in C}s_iv_i=0$ for some $C\subset[n]$ and $s\in E_C(x_1,\ldots,x_n)$, and let $H_2$ be the set of all hyperplanes of the form $\sum_{i\in C}s_iv_i=0$ for some $C\subset[n]$ and $s\not\in E_C(x_1,\ldots,x_n)$. Note that all linear equations defining the hyperplanes in $H_1\cup H_2$ have integer coefficients.

Let $V_1$ be the intersection of all hyperplanes in $H_1$, and let $V_2$ be the union of all hyperplanes in $H_2$. Since $(x_1,\ldots,x_n)\in V_1$ from definition, $V_1$ is a subspace with dimension at least 1, in which the rational points are dense. Since $(x_1,\ldots,x_n)\not\in V_2$ from definition, the subspace obtained by intersecting $V_1$ with any hyperplane in $H_2$ has a smaller dimension. As $H_2$ is finite, there is a rational point in $V_1\setminus V_2$ arbitrarily close to $(x_1,\ldots,x_n)$, and thus a point $(y_1,\ldots,y_n)$ in $V_1\setminus V_2$ with all coordinates being positive integers after clearing denominators. Therefore, $E_C(x_1,\ldots,x_n)=E_C(y_1,\ldots,y_n)$ for every $C\subset[n]$, and so $|A_>(x_1,\ldots,x_n)|=|A_>(y_1,\ldots,y_n)|$ by Corollary~\ref{cor:spitzer}.
\end{proof}

Recall that for an integer $x$, the 2-adic valuation of $x$, denoted by $\nu_2(x)$, is the largest integer $\ell$ such that $2^\ell\mid x$. We now prove an upper bound on $|E(x_1,\ldots,x_n)|$ in terms of the 2-adic valuations.
\begin{lemma}\label{lemma:localreduce12}
Let $x_1,\ldots,x_n$ be positive integers. Let $\ell=\min_{i\in[n]}\nu_2(x_i)$, let $p=|\{i\in[n]:\nu_2(x_i)=\ell\}|$, and let $q=n-p$. Let $\mathbbm{1}_p\mathbbm{2}_q$ be the vector with $p$ 1's and $q$ 2's. Then,
\[|E(x_1,\ldots,x_n)|\leq|E(\mathbbm{1}_p\mathbbm{2}_q)|.\]
\end{lemma}
\begin{proof}
Recall that $E(x_1,\ldots,x_n)$ is the number of $s\in\{-1,1\}^n$ such that $\sum_{i=1}^ns_ix_i=0$. If $p$ is odd, then both $E(x_1,\ldots,x_n)$ and $E(\mathbbm{1}_p\mathbbm{2}_q)$ are empty so we are done. 

Suppose now that $p$ is even. Let $P=\{i\in[n]:\nu_2(x_i)=\ell\}$, and let $Q=[n]\setminus P$. For every $s\in E(x_1,\ldots,x_n)$, let $I(s)=\{i\in [n]:s_i=1\}$, noting that the map $s\mapsto I(s)$ is injective. Then, $\sum_{i\in I(s)}x_i=\sum_{i\not\in I(s)}x_i=\sum_{i=1}^nx_i/2$ from definition, and it follows that the sets in $\{I(s):s\in E(x_1,\ldots,x_n)\}$ form an antichain. Furthermore, after dividing both sides of $\sum_{i\in I(s)}x_i=\sum_{i=1}^nx_i/2$ by $2^\ell$ and reducing modulo 2, we get
\[|I(s)\cap P|\equiv2^{-\ell-1}\sum_{i=1}^nx_i\pmod 2.\]
Therefore, there exists $\varepsilon\in\{0,1\}$, such that $|I(s)\cap P|\equiv\varepsilon\pmod 2$ for every $s\in E(x_1,\ldots,x_n)$.  

Let $\Pi_1$ be the poset consisting of all subsets of $P$ with size congruent to $\varepsilon$ modulo 2. Let $\Pi_2$ be the poset consisting of all subsets of $Q$. It is well known that both $\Pi_1$ and $\Pi_2$ admit symmetric chain decompositions, thus the product poset $\Pi_1\times\Pi_2$, which consists of all subsets of $[n]$ whose intersections with $P$ have size congruent to $\varepsilon$ modulo 2, has a symmetric chain decomposition as well. Therefore, $\Pi_1\times\Pi_2$ has the Sperner property, meaning that the size of every antichain within is bounded by the size of its largest rank level, where the rank of $I\in\Pi_1\times\Pi_2$ is $(|I\cap P|-\varepsilon)/2+|I\cap Q|$. The size of the level with rank $r$ is the number of $I\subset[n]$ with $|I\cap P|/2+|I\cap Q|=r+\varepsilon/2$, which is the coefficient of $z^{2r+\varepsilon}$ in $(1+z)^p(1+z^2)^q$. Since both $(1+z)^p$ and $(1+z^2)^q$ are polynomials with symmetric non-negative unimodal coefficients, their product is a polynomial with symmetric unimodal coefficients as well. Thus, $z^{p/2+q}$ has the largest coefficient in $(1+z)^p(1+z^2)^q$. 

On the other hand, from definition, $|E(\mathbbm{1}_p\mathbbm{2}_q)|$ is the constant term in $(z+z^{-1})^p(z^2+z^{-2})^{q}$. This is equal to the coefficient of $z^{p+2q}$ in $(z^2+1)^p(z^4+1)^q$, and thus equal to the coefficient of $z^{p/2+q}$ in $(1+z)^p(1+z^2)^q$. Putting everything together, $|E(x_1,\ldots,x_n)|$ is bounded from above by the size of the maximum antichain in $\Pi_1\times\Pi_2$, which is then at most the middle coefficient of $(1+z)^p(1+z^2)^q$, which is equal to $|E(\mathbbm{1}_p\mathbbm{2}_q)|$, finishing the proof.
\end{proof}

We can now reduce to the case when each of $x_1,\ldots,x_n$ is 1 or 2 with the following result. 
\begin{lemma}\label{lemma:reduce12}
Suppose for every positive integer $n$ and all non-negative integers $p,q$ with $p+q=n$, 
\[|A_>(\mathbbm{1}_p\mathbbm{2}_q)|\geq|A_>(\mathbbm{1}_n)|.\]
Then, for all positive integers $n$ and $x_1,\ldots,x_n$, 
\[|A_>(x_1,\ldots,x_n)|\geq|A_>(\mathbbm{1}_n)|=n!\binom{n-1}{\floor{(n-1)/2}}.\]
\end{lemma}
\begin{proof}
We use induction on $n$. The base cases when $n=1,2,3$ are easy to verify, so assume $n\geq4$. Let $\ell=\min_{i\in[n]}\nu_2(x_i)$, let $P=\{i\in[n]:\nu_2(x_i)=\ell\}$, and let $Q=[n]\setminus P$. Let $p=|P|$ and $q=n-p$. Note that from definition, $P\not=\varnothing$ and $Q\not=[n]$.

For every $\sigma\in\mathfrak{S}_n$, let $Q(\sigma)\subset Q$ be the union of all cycles in $\sigma$ that lie completely in $Q$. For every $U\subset Q$, let $\vec{x}_U=(x_i)_{i\in U}$, and let $\mathfrak{S}_n(U)$ be the set of permutations of $[n]\setminus U$ whose cycles all intersect $P$. Then, by using Corollary~\ref{cor:spitzer} and grouping $\sigma$ according to $Q(\sigma)$, we get
\begin{align}\label{eq:4.5}
|A_>(x_1,\ldots,x_n)|&=\sum_{\sigma\in\mathfrak{S}_n}\prod_{C\in C(\sigma)}\frac{2^{|C|}-|E_C(x_1,\ldots,x_n)|}2\nonumber\\
&=\sum_{U\subset Q}\sum_{\substack{\sigma\in\mathfrak{S}_n\\Q(\sigma)=U}}\left(\prod_{\substack{C\in C(\sigma)\\C\subset U}}\frac{2^{|C|}-|E_C(x_1,\ldots,x_n)|}2\right)\left(\prod_{\substack{C\in C(\sigma)\\C\not\subset U}}\frac{2^{|C|}-|E_C(x_1,\ldots,x_n)|}2\right)\nonumber\\
&=\sum_{U\subset Q}\left(\sum_{\sigma\in\mathfrak{S}_U}\prod_{C\in C(\sigma)}\frac{2^{|C|}-|E_C(\vec{x}_U)|}2\right)\left(\sum_{\sigma\in\mathfrak{S}_n(U)}\prod_{C\in C(\sigma)}\frac{2^{|C|}-|E_C(x_1,\ldots,x_n)|}2\right)\nonumber\\
&=\sum_{U\subset Q}|A_>(\vec{x}_U)|\left(\sum_{\sigma\in\mathfrak{S}_n(U)}\prod_{C\in C(\sigma)}\frac{2^{|C|}-|E(\vec{x}_C)|}2\right)\nonumber\\
&\geq\sum_{U\subset Q}|A_>(\mathbbm{1}_U)|\left(\sum_{\sigma\in\mathfrak{S}_n(U)}\prod_{C\in C(\sigma)}\frac{2^{|C|}-|E(\mathbbm{1}_{|C\cap P|}\mathbbm{2}_{|C\cap Q|})|}2\right),
\end{align}
where the inequality follows from induction hypothesis and Lemma~\ref{lemma:localreduce12}, and uses that $|C\cap P|\not=\varnothing$.

On the other hand, note that when $(x_1,\ldots,x_n)=\mathbbm{1}_p\mathbbm{2}_q$, the inequality above is tight as $|A_>(\vec{x}_U)|=|A_>(\mathbbm{2}_{|U|})|=|A_>(\mathbbm{1}_{|U|})|$, and $\vec{x}_C=\mathbbm{1}_{|C\cap P|}\mathbbm{2}_{|C\cap Q|}$. Thus, $|A_>(x_1,\ldots,x_n)|\geq|A_>(\mathbbm{1}_p\mathbbm{2}_q)|\geq|A_>(\mathbbm{1}_n)|$, where the last inequality follows from the assumption. 
\end{proof}

\subsection{Formula for the $1,2$-sequences}\label{sec:exact}
In this subsection, we obtain an explicit formula for $|A_>(\mathbbm{1}_p\mathbbm{2}_q)|$. To begin, we record the following result that computes the generating function for the cycle index of $\mathfrak{S}_n$. A proof can be found, for example, in~{\cite[Theorem 4.7.2]{W}}.
\begin{lemma}\label{lemma:cycleindex}
Let $z,t_1,t_2,\ldots$ be variables. Then, 
\[\sum_{n=0}^\infty\left(\sum_{\sigma\in\mathfrak{S}_n}\prod_{C\in C(\sigma)}t_{|C|}\right)\frac{z^n}{n!}=\exp\left(\sum_{k=1}^\infty\frac{t_kz^k}{k}\right).\]
\end{lemma}

For a commutative ring $R$, a formal power series $p(z)$ with coefficients in $R$, and an integer $n$, denote the coefficient of $z^n$ in $p(z)$ by $[z^n]p(z)$. Recall from the end of the proof of Lemma~\ref{lemma:localreduce12} that 
\[|E(\mathbbm{1}_p\mathbbm{2}_q)|=[z^0](z+z^{-1})^p(z^2+z^{-2})^q.\] Combining this with Corollary~\ref{cor:spitzer} and Lemma~\ref{lemma:cycleindex}, we can now compute the bivariate exponential generating function for $|A_>(\mathbbm{1}_p\mathbbm{2}_q)|$. 
\begin{lemma}\label{lemma:gen}
\[\sum_{p=0}^\infty\sum_{q=0}^\infty|A_>(\mathbbm{1}_p\mathbbm{2}_q)|\frac{u^pv^q}{p!q!}=\frac1{\sqrt{1-2u-2v}}\exp\left(\frac12[z^0]\log\left(1-u(z+z^{-1})-v(z^2+z^{-2})\right)\right).\]  
\end{lemma}
\begin{proof}
For any non-negative integers $a,b$, let 
\[w_{a,b}=\frac{2^{a+b}-|E(\mathbbm{1}_a\mathbbm{2}_b)|}2=\frac{2^{a+b}-[z^0](z+z^{-1})^a(z^2+z^{-2})^b}2.\] 
For disjoint $I,J\subset[n]$ with $I\cup J=[n]$, let $\mathbbm{1}_I\mathbbm{2}_J$ be the vector $(x_1,\ldots,x_n)$ satisfying $x_i=1$ for every $i\in I$, and $x_i=2$ for every $i\in J$. By Corollary~\ref{cor:spitzer},
\[|A_>(\mathbbm{1}_I\mathbbm{2}_J)|=\sum_{\sigma\in\mathfrak{S}_n}\prod_{C\in C(\sigma)}\frac{2^{|C|}-E_C(\mathbbm{1}_I\mathbbm{2}_J)}2=\sum_{\sigma\in\mathfrak{S}_n}\prod_{C\in C(\sigma)}w_{|C\cap I|,|C\cap J|}.\]
For every $k\geq 1$, let
\[t_k=\sum_{\substack{a,b\geq0\\a+b=k}}\binom{k}{a}w_{a,b}u^av^b.\]

Substituting the values of $t_k$ into the left hand side of Lemma~\ref{lemma:cycleindex} and setting $z=1$, we have
\begin{align*}
\sum_{n=0}^\infty\left(\sum_{\sigma\in\mathfrak{S}_n}\prod_{C\in C(\sigma)}t_{|C|}\right)\frac1{n!}&=\sum_{n=0}^\infty\left(\sum_{\sigma\in\mathfrak{S}_n}\prod_{C\in C(\sigma)}\left(\sum_{\substack{a_C,b_C\geq0\\a_C+b_C=|C|}}\binom{|C|}{a_C}w_{a_C,b_C}u^{a_C}v^{b_C}\right)\right)\frac1{n!}\\
&=\sum_{p=0}^\infty\sum_{q=0}^\infty\left(\sum_{\sigma\in\mathfrak{S}_{p+q}}\sum_{\substack{(a_C),(b_C)\\a_C+b_C=|C|\\\sum a_C=p,\sum b_C=q}}\prod_{C\in C(\sigma)}\binom{|C|}{a_C}w_{a_C,b_C}\right)\frac{u^pv^q}{(p+q)!}\\
&=\sum_{p=0}^\infty\sum_{q=0}^\infty\sum_{I\in\binom{[n]}{p}}|A_>(\mathbbm{1}_I\mathbbm{2}_{[n]\setminus I})|\frac{u^pv^q}{(p+q)!}\\
&=\sum_{p=0}^\infty\sum_{q=0}^\infty\binom{p+q}{p}|A_>(\mathbbm{1}_p\mathbbm{2}_q)|\frac{u^pv^q}{(p+q)!}=\sum_{p=0}^\infty\sum_{q=0}^\infty|A_>(\mathbbm{1}_p\mathbbm{2}_q)|\frac{u^pv^q}{p!q!},
\end{align*}
where the fourth sum after the second equality sign is summing over all sequences of non-negative integers $(a_C)$ and $(b_C)$ indexed by $C\in C(\sigma)$, such that $a_C+b_C=|C|$ for every $C\in C(\sigma)$, $\sum_{C\in C(\sigma)}a_C=p$, and $\sum_{C\in C(\sigma)}b_C=q$.

On the other hand, substituting the values of $t_k$ into the right hand side of Lemma~\ref{lemma:cycleindex} and again setting $z=1$, we get
\begin{align*}
\exp\left(\sum_{k=1}^\infty\frac{t_k}{k}\right)&=\exp\left(\sum_{k=1}^\infty\frac1{k}\sum_{\substack{a,b\geq0\\a+b=k}}\binom{k}{a}w_{a,b}u^av^b\right)=\exp\left(\sum_{\substack{a,b\geq0\\a+b\geq1}}\frac{\binom{a+b}a}{a+b}w_{a,b}u^av^b\right)\\
&=\exp\left(\sum_{\substack{a,b\geq0\\a+b\geq1}}\frac{\binom{a+b}a}{a+b}\left(\frac{2^{a+b}-[z^0](z+z^{-1})^a(z^2+z^{-2})^b}2\right)u^av^b\right).
\end{align*}
Recall that
\[-\log(1-x-y)=\sum_{\substack{a,b\geq0\\a+b\geq1}}\frac{\binom{a+b}a}{a+b}x^ay^b.\]
By substituting first $x=2u$, $y=2v$, and then $x=(z+z^{-1})u$, $y=(z^2+z^{-2})v$, we get
\[\sum_{\substack{a,b\geq0\\a+b\geq1}}\frac{\binom{a+b}a}{a+b}2^{a+b}u^av^b=-\log(1-2u-2v),\]
\[\sum_{\substack{a,b\geq0\\a+b\geq1}}\frac{\binom{a+b}a}{a+b}(z+z^{-1})^a(z^2+z^{-2})^bu^av^b=-\log\left(1-(z+z^{-1})u-(z^2+z^{-2})v\right).\]
Therefore,
\[\exp\left(\sum_{k=1}^\infty\frac{t_k}{k}\right)=\frac1{\sqrt{1-2u-2v}}\exp\left(\frac12[z^0]\log\left(1-(z+z^{-1})u-(z^2+z^{-2})v\right)\right),\]
which completes the proof.
\end{proof}

To finish finding a formula for $|A_>(\mathbbm{1}_p\mathbbm{2}_q)|$, we need to extract the coefficient of $u^pv^q$ from the generating function obtained in Lemma~\ref{lemma:gen}. We first substitute $u$ and $v$ with carefully chosen new formal power series in order to remove the cumbersome coefficient extraction of $z^0$ from inside the exponent, then use the following Lagrange-Good inversion formula to express the desired coefficient of $u^pv^q$ in terms of $r^ps^q$. For simplicity, this is only stated in the two-variable version that we need.
\begin{lemma}[Lagrange-Good inversion formula~{\cite[Theorem 1]{Ge}}]\label{lemma:invert}
Let $\phi_1(r,s)$ and $\phi_2(r,s)$ be formal power series with non-zero constant terms. Define formal power series $u(r,s)$ and $v(r,s)$ by $u(r,s)=r/\phi_1(r,s)$ and $v(r,s)=s/\phi_2(r,s)$. Then, for any formal power series $\phi(u,v)$ and any $p,q\geq0$,
\[[u^pv^q]\phi(u,v)=[r^ps^q]\phi(u(r,s),v(r,s))\phi_1(r,s)^p\phi_2(r,s)^q\det\begin{pmatrix}
    1-\dfrac{r}{\phi_1}\dfrac{\partial\phi_1}{\partial r} & -\dfrac{s}{\phi_1}\dfrac{\partial\phi_1}{\partial s}\\-\dfrac{r}{\phi_2}\dfrac{\partial\phi_2}{\partial r} & 1-\dfrac{s}{\phi_2}\dfrac{\partial\phi_2}{\partial s}
\end{pmatrix}.\]
\end{lemma}

\begin{lemma}\label{lemma:extract}
For all non-negative integers $p,q$ with $p+q\geq 1$,
\begin{align*}
|A_>(\mathbbm{1}_p\mathbbm{2}_q)|&=p!q![r^ps^q](1+r^2+s^2)^{p+q-1}\frac{(1+s)(1+r-s)}{(1-s)^{p+1}}\\
&=p!q!\binom{p+q-1}{\floor{p/2}}[s^q]\frac{(1+s^2)^{\ceil{p/2}+q-1}(1+s)}{(1-s)^{2\ceil{p/2}}}
\end{align*}
\end{lemma}
\begin{proof}
Let $\phi_1(r,s)=(1+r^2+s^2)/(1-s)$ and $\phi_2(r,s)=1+r^2+s^2$. Let $u=r/\phi_1(r,s)$ and $v=s/\phi_2(r,s)$. Then, the assumptions of Lemma~\ref{lemma:invert} are satisfied. 

Let $\phi(u,v)$ be the formal power series on the right hand side of Lemma~\ref{lemma:gen}. We compute $\phi(u(r,s),v(r,s))$ as follows.
\[\frac1{\sqrt{1-2u(r,s)-2v(r,s)}}=\frac1{\sqrt{1-\frac{2r(1-s)+2s}{1+r^2+s^2}}}=\frac{\sqrt{1+r^2+s^2}}{1-r-s},\]
\begin{align*}
[z^0]\log\left(1-u(r,s)(z+z^{-1})-v(r,s)(z^2+z^{-2})\right)&=[z^0]\log\left(\frac{(1-rz-sz^2)(1-rz^{-1}-sz^{-2})}{1+r^2+s^2}\right)\\
&=-\log(1+r^2+s^2).
\end{align*}
where the last equation follows as both $\log(1-rz-sz^2)$ and $\log(1-rz^{-1}-sz^{-2})$ have constant term zero. Therefore, 
\[\phi(u(r,s),v(r,s))=\frac{\sqrt{1+r^2+s^2}}{1-r-s}\exp\left(-\frac12\log(1+r^2+s^2)\right)=\frac1{1-r-s}.\]

Now, by Lemma~\ref{lemma:invert}, we get
\begin{align*}
|A_>(\mathbbm{1}_p\mathbbm{2}_q)|&=p!q![u^pv^q]\phi(u,v)\\
&=p!q![r^ps^q]\phi(u(r,s),v(r,s))\phi_1(r,s)^p\phi_2(r,s)^q\det\begin{pmatrix}
    1-\dfrac{r}{\phi_1}\dfrac{\partial\phi_1}{\partial r} & -\dfrac{s}{\phi_1}\dfrac{\partial\phi_1}{\partial s}\\-\dfrac{r}{\phi_2}\dfrac{\partial\phi_2}{\partial r} & 1-\dfrac{s}{\phi_2}\dfrac{\partial\phi_2}{\partial s}
\end{pmatrix}\\
&=p!q![r^ps^q]\frac1{1-r-s}\frac{(1+r^2+s^2)^p}{(1-s)^p}(1+r^2+s^2)^q\frac{(1+s)(1-r-s)(1+r-s)}{(1-s)(1+r+s)}\\
&=p!q![r^ps^q](1+r^2+s^2)^{p+q-1}\frac{(1+s)(1+r-s)}{(1-s)^{p+1}}.
\end{align*}

Finally, note that if $p$ is even, then the $r^p$ term must come from $(1+r^2+s^2)^{p+q-1}$ by selecting the $r^2$ term $p/2$ times, while if $p$ is odd, then the $r^p$ term must be formed by selecting the $r^2$ term $(p-1)/2$ times from $(1+r^2+s^2)^{p+q-1}$, and selecting the $r$ term from $1+r-s$. Using floor and ceiling functions appropriately then gives the unifying formula as stated. 
\end{proof}

\subsection{Optimization}\label{sec:opt}
In this subsection, we use the formula for $|A_>(\mathbbm{1}_p\mathbbm{2}_q)|$ obtained in Lemma~\ref{lemma:extract} to determine when $|A_>(\mathbbm{1}_p\mathbbm{2}_q)|$ is minimized subject to $p+q=n$. 

\begin{lemma}\label{lemma:optimize}
Let $n\geq1$ and $p,q\geq0$ be integers with $n=p+q$. Then, $|A_>(\mathbbm{1}_p\mathbbm{2}_q)|\geq n!\binom{n-1}{\floor{(n-1)/2}}$. Moreover, equality holds only when $(p,q)=(n,0)$ and $(0,n)$, and when $(p,q)=(n-1,1)$ if $n$ is odd.  
\end{lemma}
\begin{proof} 
For all non-negative integers $c,d$, let 
\[B(c,d)=[s^d]\frac{(1+s^2)^{c+d-1}(1+s)}{(1-s)^{2c}},\]
so that by Lemma~\ref{lemma:extract}, $|A_>(\mathbbm{1}_p\mathbbm{2}_q)|=p!q!\binom{n-1}{\floor{p/2}}B(\ceil{p/2},q)$. For the four possible parity combinations of $p,q$, proving $|A_>(\mathbbm{1}_p\mathbbm{2}_q)|\geq n!\binom{n-1}{\floor{(n-1)/2}}$ is equivalent to proving the following for all relevant $c,d$.
\begin{alignat*}{2}
B(c,2d)\geq E_1(c,2d)&:=\frac{\binom{2c+2d}{2c}\binom{2c+2d-1}{c+d-1}}{\binom {2c+2d-1}c}, &&p=2c, q=2d.\\
B(c,2d)\geq E_2(c,2d)&:=\frac{\binom{2c+2d-1}{2c-1}\binom{2c+2d-2}{c+d-1}}{\binom{2c+2d-2}{c-1}}, &&p=2c-1, q=2d.\\
B(c,2d-1)\geq O_1(c,2d-1)&:=\frac{\binom{2c+2d-1}{2c}\binom{2c+2d-2}{c+d-1}}{\binom {2c+2d-2}c}, \quad\quad&&p=2c, q=2d-1.\\
B(c,2d-1)\geq O_2(c,2d-1)&:=\frac{\binom{2c+2d-2}{2c-1}\binom{2c+2d-3}{c+d-2}}{\binom{2c+2d-3}{c-1}}, &&p=2c-1, q=2d-1.
\end{alignat*}
Straightforward computations show that $O_2(c,2d-1)=\frac{2c+2d-2}{2c+2d-1}O_1(c,2d-1)<O_1(c,2d-1)$, while $E_1(c,2d)=E_2(c,2d)$. Let $O(c,2d-1)=O_1(c,2d-1)$ for all $c\geq0$ and $d\geq 1$. Let $E(0,2d)=E_1(0,2d)$ for all $d\geq0$, and let 
\begin{equation}\label{eq:1}
E(c,2d)=\frac12(O(c,2d+1)-O(c-1,2d+1))
\end{equation}
for all $c\geq 1$ and $d\geq0$. Then, a direct computation gives for all $c\geq 1$ and $d\geq 0$,
\begin{align*}
\frac{E(c,2d)}{E_1(c,2d)}&=\frac{O(c,2d+1)}{2E_1(c,2d)}-\frac{O(c-1,2d+1)}{2E_1(c,2d)}\\
&=\frac{(2c+2d+1)(c+2d)}{2(c+d)(2d+1)}-\frac{2c-1}{2(2d+1)}=1+\frac{d}{2(c+d)(2d+1)}.
\end{align*}
Similarly, after some tedious calculations, we have that for any $c\geq2$ and $d\geq1$,
\begin{align}\label{eq:2}
&\phantom{{}={}}\frac{E(c-1,2d)+2O(c,2d-1)-E(c,2d)}{E_1(c,2d)}\nonumber\\
&=\frac{(2c-1)(4cd+4d^2+2c-d-2)}{2(2d+1)(c+2d-1)(2c+2d-1)}+\frac{2d}{c+2d-1}-1-\frac{d}{2(c+d)(2d+1)}\nonumber\\
&=\frac{d(4c+5d-1)}{2(c+d)(2d+1)(c+2d-1)(2c+2d-1)}>0.
\end{align}

With these in mind, we now use induction on $c+2d$ to prove that $B(c,2d-1)\geq O(c,2d-1)$ for all $c\geq0$ and $d\geq1$, with the inequality being strict when $c\geq 1$ and $d\geq2$, and that $B(c,2d)\geq E(c,2d)$ for all $c\geq0$ and $d\geq0$, with the inequality being strict when $c\geq 1$ and $d\geq 1$. The cases when $c=0$, $c=1$, or $d=0$ can be verified by computing the relevant values directly.

Note that for $c,d\geq1$, we have the useful recurrence
\begin{align*}
B(c-1,d)+2B(c,d-1)&=[s^d]\left(\frac{(1+s^2)^{c+d-2}(1+s)}{(1-s)^{2c-2}}+\frac{2(1+s^2)^{c+d-2}(1+s)s}{(1-s)^{2c}}\right)\\
&=[s^d]\frac{(1+s^2)^{c+d-2}(1+s)((1-s)^2+2s)}{(1-s)^{2c}}\\
&=\frac{(1+s^2)^{c+d-1}(1+s)}{(1-s)^{2c}}=B(c,d).
\end{align*}
Using this and the induction hypothesis, we have in the odd case for all $c\geq2$ and $d\geq1$,
\[B(c,2d-1)=B(c-1,2d-1)+2B(c,2d-2)\geq O(c-1,2d-1)+2E(c,2d-2)=O(c,2d-1),\]
where the last equality uses~(\ref{eq:1}). Moreover, $B(c-1,2d-1)>O(c-1,2d-1)$ when $c,d\geq2$ from induction hypothesis, so $B(c,2d-1)>O(c,2d-1)$ as well.
Similarly, in the even case, using~(\ref{eq:2}), for all $c\geq 2$ and $d\geq 1$, 
\[B(c,2d)=B(c-1,2d)+2B(c,2d-1)\geq E(c-1,2d)+2O(c,2d-1)>E(c,2d).\]

Recall that $O(c,2d-1)=O_1(c,2d-1)>O_2(c,2d-1)$, and $E(c,2d)>E_1(c,2d)=E_2(c,2d)$ unless $c=0$ or $d=0$. Combined with the beginning of the proof, we have thus proved that $|A_>(\mathbbm{1}_p\mathbbm{2}_q)|\geq n!\binom{n-1}{\floor{(n-1)/2}}$. Furthermore, if $q$ is odd and equality holds, then we need $B(c,2d-1)=O(c,2d-1)=O_1(c,2d-1)>O_2(c,2d-1)$, which implies that either $p=c=0$ and $(p,q)=(0,n)$, or $q=d=1$ and $p$ is even, and so $(p,q)=(n-1,1)$ with $n$ odd. If $q$ is even and equality holds, then we need $B(c,2d)=E(c,2d)=E_1(c,2d)=E_2(c,2d)$, so either $p=c=0$ and $(p,q)=(0,n)$, or $q=d=0$ and $(p,q)=(n,0)$. Finally, equality does hold in each of these cases by Lemma~\ref{lemma:tight}, which finishes the proof.
\end{proof}

We are now ready to put everything together and prove Theorem~\ref{thm:mainstrong}.
\begin{proof}[Proof of Theorem~\ref{thm:mainstrong}]
Let $x_1,\ldots,x_n>0$. By Lemma~\ref{lemma:reduceint}, there exist positive integers $y_1,\ldots,y_n$ such that $E_C(x_1,\ldots,x_n)=E_C(y_1,\ldots,y_n)$ for every $C\subset[n]$, and so $|A_>(x_1,\ldots,x_n)|=|A_>(y_1,\ldots,y_n)|$. By Lemma~\ref{lemma:reduce12} and Lemma~\ref{lemma:optimize}, $|A_>(y_1,\ldots,y_n)|\geq n!\binom{n-1}{\floor{(n-1)/2}}$, so $|A_>(x_1,\ldots,x_n)|\geq n!\binom{n-1}{\floor{(n-1)/2}}$ as well, proving the lower bound.

Now suppose $|A_>(x_1,\ldots,x_n)|=n!\binom{n-1}{\floor{(n-1)/2}}$, then $|A_>(y_1,\ldots,y_n)|=n!\binom{n-1}{\floor{(n-1)/2}}$ as well. From the proof of Lemma~\ref{lemma:reduce12}, if $\ell=\min_{i\in[n]}\nu_2(y_i)$, $p=|\{i\in[n]:\nu_2(y_i)=\ell\}|$, and $q=n-p$, then $n!\binom{n-1}{\floor{(n-1)/2}}=|A_>(y_1,\ldots,y_n)|\geq|A_>(\mathbbm{1}_p\mathbbm{2}_q)|$, so $|A_>(\mathbbm{1}_p\mathbbm{2}_q)|=n!\binom{n-1}{\floor{(n-1)/2}}$ as well. By Lemma~\ref{lemma:optimize}, and since $p\geq 1$, we must have $(p,q)=(n,0)$, or $n$ is odd and $(p,q)=(n-1,1)$. 

If $(p,q)=(n,0)$, then we claim that $x_1=\cdots=x_n$. Indeed, since $|A_>(y_1,\ldots,y_n)|=|A_>(\mathbbm{1}_p\mathbbm{2}_q)|$, the inequality~(\ref{eq:4.5}) in the proof of Lemma~\ref{lemma:reduce12} must actually be an equality. In this case, we have $P=[n]$ and $Q=\varnothing$, so for every $C\subset [n]$, $|E(\vec{y}_C)|=|E(\mathbbm{1}_{|C|})|$. In particular, for any distinct $i,j\in[n]$, $|E(y_i,y_j)|=|E(1,1)|=2$, so $y_i=y_j$, as required. Since $E_C(x_1,\ldots,x_n)=E_C(y_1,\ldots,y_n)$ for every $C\subset[n]$, we have $x_i=x_j$ for all distinct $i,j\in[n]$, again by taking $C=\{i,j\}$.  

If $n$ is odd and $(p,q)=(n-1,1)$, without loss of generality assume that $\nu_2(y_n)$ is the unique larger one, so that $P=[n-1]$ and $Q=\{n\}$. Again, the inequality~(\ref{eq:4.5}) in the proof of Lemma~\ref{lemma:reduce12} must actually be an equality. By taking $\sigma=(i j)$ for all distinct $i,j\in[n-1]$, we get $y_i=y_j$. By taking $\sigma=(i j n)$ for all distinct $i,j\in[n-1]$, we get that $|E(y_i,y_j,y_n)|=|E(1,1,2)|=2$, so $y_n=2y_i=2y_j$. Thus, $y_1=\cdots=y_{n-1}=y_n/2$. Since $E_C(x_1,\ldots,x_n)=E_C(y_1,\ldots,y_n)$ for every $C\subset[n]$, we also have $x_1=\cdots=x_{n-1}=x_n/2$, which finishes the proof.  
\end{proof}

\section{Proof of Corollaries}\label{sec:cor}
In this section, we prove Corollary~\ref{cor:pmain} and Corollary~\ref{cor:piid}.
\begin{proof}[Proof of Corollary~\ref{cor:pmain}]
Let $M$ be the random variable counting the number of $k\in[n]$ such that $X_k\not=0$. Conditioning on $M=m$, let $X_1',\ldots,X_m'$ be the sequence of variables among $X_1,\ldots,X_n$ which are non-zero and listed in the same order. Let $S_k'=\sum_{i=1}^kX_i'$ for every $k\in[m]$. Since $(X_1,\ldots,X_n)$ is exchangeable and sign-invariant, so is $(X_1',\ldots,X_m')$. Thus, by Theorem~\ref{thm:esatom} and Theorem~\ref{thm:main1}, we have
\[\frac{1}{2^m}\binom{m-1}{\floor{(m-1)/2}}\leq\mathbb{P}(S_1',\ldots,S_m'>0)\leq\frac{1}{4^m}\binom{2m}{m}\leq\mathbb{P}(S_1',\ldots,S_m'\geq0)\leq\frac{1}{2^m}\binom{m}{\floor{m/2}}.\]

For weak persistence, after conditioning on $M=m$, we have $S_1,\ldots,S_n\geq0$ if and only if $S_1',\ldots,S_m'\geq0$. Thus,
\begin{align*}
\mathbb{P}(S_1,\ldots,S_n\geq0)&=\sum_{m=0}^n\mathbb{P}(M=m)\mathbb{P}(S_1',\ldots,S_m'\geq0)\\
&\geq\sum_{m=0}^n\mathbb{P}(M=m)\frac{1}{4^m}\binom{2m}{m}=\mathbb{E}\left(\frac{1}{4^M}\binom{2M}{M}\right),  
\end{align*}
and similarly
\[\mathbb{P}(S_1,\ldots,S_n\geq0)\leq\sum_{m=0}^n\mathbb{P}(M=m)\frac{1}{2^m}\binom{m}{\floor{m/2}}=\mathbb{E}\left(\frac{1}{2^M}\binom{M}{\floor{M/2}}\right).\]

For strong persistence, note that conditioning on $M=m$, $S_1,\ldots,S_n>0$ if and only if $X_1\not=0$ and $S_1',\ldots,S_m'>0$. Furthermore, since $(X_1,\ldots,X_n)$ is exchangeable and sign-invariant, $\mathbb{P}(X_1\not=0\mid M=m)=m/n$, and $\mathbb{P}(S_1',\ldots,S_m'>0)=\mathbb{P}(S_1',\ldots,S_m'>0\mid X_1\not=0)$. Thus, we have
\begin{align*}
\mathbb{P}(S_1,\ldots,S_n>0)&=\sum_{m=1}^n\mathbb{P}(M=m)\mathbb{P}(X_1\not=0\mid M=m)\mathbb{P}(S_1,\ldots,S_n>0\mid M=m,X_1\not=0)\\
&=\sum_{m=1}^n\mathbb{P}(M=m)\cdot\frac mn\cdot\mathbb{P}(S_1',\ldots,S_m'>0)\\
&\geq\sum_{m=1}^n\mathbb{P}(M=m)\cdot\frac mn\cdot\frac{1}{2^m}\binom{m-1}{\floor{(m-1)/2}}=\frac1n\mathbb{E}\left(\frac M{2^M}\binom{M-1}{\floor{(M-1)/2}}\right),
\end{align*}
where we interpret $\frac{M}{2^M}\binom{M-1}{\floor{(M-1)/2}}$ as 0 when $M$ takes the value 0. Similarly
\begin{align*}
\mathbb{P}(S_1,\ldots,S_n>0)&=\sum_{m=1}^n\mathbb{P}(M=m)\cdot\frac mn\cdot\mathbb{P}(S_1',\ldots,S_m'>0)\\
&\leq\sum_{m=1}^n\mathbb{P}(M=m)\cdot\frac mn\cdot\frac{1}{4^m}\binom{2m}{m}=\frac1n\mathbb{E}\left(\frac{M}{4^M}\binom{2M}{M}\right).
\end{align*}

We now show that for every given random variable $M$, these bounds are all tight. Let $Y_1$ be the random variable taking values $-1$ and 1 with probability $1/2$ each, and let $Y_2$ be any symmetric atomless random variable. First, define $(X_1,\ldots,X_n)$ as follows. Conditioning on $M=m$, pick a uniformly random $I\subset[n]$ of size $m$, let $(X_i)_{i\in I}$ be independent identically distributed copies of $Y_1$, and let $X_i$ be 0 for every $i\not\in I$. Then, the equality case of Theorem~\ref{thm:main1} shows that $(X_1,\ldots,X_n)$ attains the weak persistence upper bound and strong persistence lower bound. For the weak persistence lower bound and strong persistence upper bound, the same definition except with $Y_2$ in place of $Y_1$ achieves both using the equality case of Theorem~\ref{thm:esatom}.
\end{proof}

\begin{proof}[Proof of Corollary~\ref{cor:piid}]
The case when $p=0$ is trivial, so assume $p>0$. Let $X'$ be the symmetric random variable obtained from $X$ by conditioning on $X\not=0$. More precisely, for every $x\in\mathbb{R}$, let
\[\mathbb{P}(X'\leq x)=\mathbb{P}(X\leq x\mid X\not=0).\]
Let $X_1',\ldots,X_m'$ be independent identically distributed copies of $X'$, and define $S_k'=\sum_{i=1}^kX_i'$ for every $k\in[m]$. Note that $(X_1',\ldots,X_m')$ is exchangeable and sign-invariant, so by Theorem~\ref{thm:esatom} and Theorem~\ref{thm:main1}, we again have
\begin{equation}\label{eq:piid}
\frac{1}{2^m}\binom{m-1}{\floor{(m-1)/2}}\leq\mathbb{P}(S_1',\ldots,S_m'>0)\leq\frac1{4^m}\binom{2m}m\leq\mathbb{P}(S_1',\ldots,S_m'\geq0)\leq\frac{1}{2^m}\binom{m}{\floor{m/2}}.
\end{equation}

Let $M$ be the random variable counting the number of $k\in[n]$ with $X_k\not=0$. Since $M\sim\text{Bin}(n,2p)$, $\mathbb{P}(M=m)=\binom nm(2p)^m(1-2p)^{n-m}$. Like in the proof of Corollary~\ref{cor:pmain}, for weak persistence we have
\[\mathbb{P}(S_1,\ldots,S_n\geq0)=\sum_{m=0}^n\binom nm(2p)^m(1-2p)^{n-m}\mathbb{P}(S_1',\ldots,S_m'\geq0),\]
and for strong persistence we have
\begin{align*}
\mathbb{P}(S_1,\ldots,S_n>0)&=\sum_{m=1}^n\mathbb{P}(M=m)\mathbb{P}(X_1\not=0\mid M=m)\mathbb{P}(S_1,\ldots,S_n>0\mid M=m,X_1\not=0)\\
&=\sum_{m=1}^n\binom{n}{m}(2p)^m(1-2p)^{n-m}\frac{m}{n}\mathbb{P}(S_1,\ldots,S_n>0\mid M=m,X_1\not=0)\\
&=\sum_{m=1}^n\binom{n-1}{m-1}(2p)^m(1-2p)^{n-m}\mathbb{P}(S_1',\ldots,S_m'>0).
\end{align*}
Combining these with (\ref{eq:piid}), all four bounds are proved. 

Let $Y_1$ and $Y_2$ be the random variables defined at the end of the proof of Corollary~\ref{cor:pmain}. Using the equality case of Theorem~\ref{thm:esatom} and Theorem~\ref{thm:main1} again, we see that if $X$ is defined so that $X'=Y_1$, then the weak persistence upper bound and strong persistence lower bound are attained, while if $X$ is defined so that $X'=Y_2$, then the weak persistence lower bound and strong persistence upper bound are attained. 
\end{proof}

\section*{Acknowledgements}
The main results presented in this paper are the product of a human-AI collaborative research process utilizing the Google DeepMind AI co-mathematician \cite{zheng2026ai} (for Theorem~\ref{thm:mainweak}) and ChatGPT (for Theorem~\ref{thm:mainstrong}). The writing and presentation of this paper are the work of the authors, who take full responsibility for the content.

\bibliographystyle{plain}
\bibliography{bibliography}
\end{document}